\documentclass[12pt]{amsart}
\usepackage{amsfonts}
\usepackage{amsmath}
\usepackage{amsthm}
\usepackage{amssymb}
\usepackage{amscd}
\usepackage{verbatim}
\usepackage{color}
\usepackage{url}

\usepackage{tikz}

\usepackage{graphicx}
\usepackage{microtype}
\usepackage{enumerate}
\usepackage{pinlabel}

\theoremstyle{definition}

\newtheorem{thm}{Theorem}[section]
\newtheorem{cor}[thm]{Corollary}
\newtheorem{prop}[thm]{Proposition}
\newtheorem{lem}[thm]{Lemma}
\newtheorem{defi}[thm]{Definition}

\newtheorem{exam}[thm]{Example}

\newcommand{\R}{\mathbb{R}}
\newcommand{\Z}{\mathbb{Z}}
\newcommand{\N}{\mathbb{N}}

\newcommand{\vs}{\vspace{0,3cm}}
\newcommand{\vsp}{\vspace{0,15cm}}

\title[The space of relative orders and  Morris indicability theorem]{The space of relative orders and a generalization of Morris indicability theorem}
\author{Yago Antol\'in} 
\author{Crist\'obal Rivas}
\address{Dpto. de Matem\'aticas, Universidad Aut\'onoma de Madrid and Instituto de
Ciencias Matemáticas, CSIC-UAM-UC3M-UCM.}
\address{Dpto. Matem\'atica y C.C. Universidad de Santiago de Chile, Alameda 3363, Estaci\'on Central, Santiago, Chile}
\date{}
\begin{document}

\thanks{The  first author is supported by the Juan de la Cierva grant IJCI-2014-22425, and acknowledges partial
support from the Spanish Government through grant number MTM2014-54896 and through the “Severo Ochoa Programme for Centres of Excellence in R\&{}D” (SEV-2015‐0554). The second author is partially supported by FONDECYT 1181548}

\begin{abstract}
We introduce the space of relative orders on a group and show that it is compact whenever the group is finitely generated. 
We use this to show that if $G$ is a finitely generated group acting by order preserving homeomorphism of on the line, then if  some stabilizer of a point is proper and co-amenable subgroup, then $G$ surjects onto $\mathbb{Z}$. This is a generalization of a theorem of Morris.

\end{abstract}


\maketitle

\section{Introduction}

In \cite{Mo}, Dave  Morris proved the following  celebrated theorem

\begin{thm}[Morris] \label{teo Morris} {\em  Let $G$ be a finitely generated and amenable group acting non trivially by order preserving homeomorphisms on the line. Then $G$ surjects onto $\Z$.}
\end{thm}

His proof exploits a strong connection between faithful group actions on the line and left-multiplication invariant total orders ({\em left-orders} for short) on groups (see \cite{GOD} for a general introduction on this relationship). Crucial to his proof, is the fact that $\mathcal{LO}(G)$, the space of all left-orders on a group $G$, is compact \cite{sikora} and there is a natural $G$-action on it by homeomorphisms. Thus there is a $G$-invariant probability measure on $\mathcal{LO}(G)$ whenever $G$ is amenable. Finally, Morris shows that {\em almost every} order in the support of the invariant probability is a  left-order of Conradian type\footnote{A left-order $\preceq$ is of Conradian type if whenever $id\preceq f \prec g$ then $g\prec fg^n$ for some $n\geq1$. See \cite{conrad}.}, which implies that $G$ surjects onto $\Z$ whenever $G$ is finitely generated \cite{conrad}.

In this note we extend Morris's result to groups acting on the line that are not necessarily amenable but contain large point-stabilizers. More precisely we show

\begin{thm} \label{teo main} {\em Let $G$ be a finitely generated group acting on the line by order preserving homeomorphisms. Suppose there is a point $p\in \R$ whose stabilizer in $G$ is a proper and co-amenable subgroup. Then $G$ surjects onto $\Z$.}
\end{thm}

Recall that $H\leq G$ is {\bf co-amenable} if whenever $G$ acts by homeomorphisms of a compact space in such a way that $H$ preserves a probability measure, then $G$ also preserves a probability measure (as with amenability, there are many equivalent definitions of co-amenability, see \cite{eymard, monod popa}). So, amenable groups are precisely groups on which the identity is co-amenable. In particular, Morris' theorem follows from ours as the case that the kernel of the action is co-amenable.

Given $G$ a finitely generated group acting on the line by order preserving  homeomorphisms and $H=Stab_G(p)$ a proper subgroup, a natural reaction is to try to show  that $H$ preserves a probability measure when acting on $\mathcal{LO}(G)$. This is however, not true in general as we explain in Example \ref{ejemplo}.  In fact, the essence of our work is to show that $\mathcal{LO}(G)$ sits inside a larger space, which we call the {\em space of relative orders}, which is also a compact space enjoying a $G$-action by homeomorphisms and, moreover, it is enlarged enough so that the $H$-action on it will have fixed points.

\begin{defi}[Folklore] {\em A  proper subgroup $C$ of a group $G$ is called {\em relatively convex} if  there is a total order $\preceq$ on $G/C$ which is invariant under left-multiplication by $G$, meaning that if $fC\preceq gC$, then $hfC\preceq hgC$ for all $f,g,h$ in $G$. (In particular, $C$ has infinite index in $G$.)

We say that $\preceq$ is a relative order of $G$ (with respect  to $C$) and denote by $\mathcal{O}_{Rel}(G)$ the set of all relative orders of $G$.}
\end{defi}

 In Section \ref{sec relative orders} we will show

\begin{thm} \label{teo Rel}{\em There is a natural topology on $\mathcal{O}_{Rel}(G)$ that makes it compact whenever $G$ is finitely generated. Moreover, $G$ acts naturally on $\mathcal{O}_{Rel}(G)$ by homeomorphisms and, if $C$ is a relatively convex subgroup of $G$, then $C$ fixes a point in $\mathcal{O}_{Rel}(G)$.}
\end{thm}

Here is one easy recipe to build relative orders. Start with a non-trivial action of a group $G$ on the line and take $x\in \R$ a point which is not globally fixed by $G$. Then, we can declare $f\preceq g$ if and only if $g(x)\leq f(x)$. This is a relative order of $G$ whose corresponding relative convex subgroup is $Stab_G(x)$. In Section \ref{sec relative orders}, we will also show that any relative order on $G$ can be obtained from the above procedure. Thus, if $G$ is finitely generated, $\mathcal{O}_{Rel}(G)$ can be regarded as a compactification of the space of non-trivial actions of $G$ on the line in the same way as $\mathcal{LO}(G)$ can be regarded as a compactification of the space of faithful actions of $G$ on the line.

There is just one small caveat in order to implement Morris' strategy with $\mathcal{O}_{Rel}(G)$ in place of $\mathcal{LO}(G)$. 
Morris in his paper, after finding an order that is  Conradian, gives an alternative argument for the indicability that avoids Conrad's theorem.
In our case, we are forced to avoid Conrad's theorem since  there is no available analog in the context of  relative orders on groups.
To overcome this issue we use the concept of {\em crossings} for a group acting on a totally ordered space (see Definition \ref{def crossing}). This concept was introduced by Beklaryan \cite{beklaryan} in the case of group acting on the line, but it was Navas \cite{navas} and later Navas and the second author \cite{Na-Ri} who realized that the absence of crossings in the left multiplication action of $G$ on the totally ordered space $(G,\preceq)$ is a characterization of the Conradian property for the left-order $\preceq$. In Section \ref{sec coamenable} we will generalize this observation by showing that actions without crossings of a finitely generated group on a totally ordered space always entails a surjection onto $\Z$ (Proposition \ref{prop crossings}). In addition to that, we will show that  {\em almost every} relative order in the support of a $G$-invariant probability measure on $\mathcal{O}_{Rel}(G)$ corresponds to an action of $G$ on $(G/C,\preceq)$ without crossings, thus showing Theorem \ref{teo main}.
 



\section{The space of relative orders}
\label{sec relative orders}


We will use the following characterization of relative orders, see \cite[Corollary 5.1.5]{KopytovMedvedev} or \cite{ADS}. For completeness we provide a proof.

\begin{lem}\label{lem:charconvex}
{\em The subgroup $C$ is relatively convex in $G$ if and only if there is a  semigroup $P\neq\{id\}$ such that \begin{enumerate}
\item[$i)$] $CPC\subset P$ and 
\item[$ii)$] $G=P\sqcup P^{-1} \sqcup C$ (disjoint union). 
\end{enumerate}
If $P$ is a semigroup satisfying $i)$ and $ii)$ for some subgroup $C$, we say that  $P$ a {\em relative cone} with respect to $C$.}

\end{lem}

\proof  Let $C$ be a relatively convex subgroup of $G$ and denote by $\preceq$ a $G$-invariant total order of $G/C$. Then $P=\{g\in G\mid gC\succ C\}$ is a relative cone. Conversely if $P$ is a relative cone respect to $C$, define $f C\prec gC$ if and only if $f^{-1}g\in P$. Observe that this definition is independent of the coset representative and that two cosets $fC$ and $gC$ are not comparable under $\prec$, if and only if $f^{-1}g\in C$. $\hfill\square$

\vsp

Some authors, {\em e.g.} \cite{KopytovMedvedev}, allow a relatively convex subgroup to be the whole ambient group. The reason why we restrict our attention only to proper subgroups, is because under that restriction relative orders correspond to non-trivial actions of the group on the line,  in the same way as total left multiplication invariant orders (left orders for short) correspond to faithful actions on the line by orientation-preserving homeomorphisms, see \cite{ghys}. More precisely we have

\begin{prop} \label{prop real dyn} {\em Suppose that $G$ is a countable group. Then $G$ admits a relative order if and only if it admits a non-trivial action by orientation-preserving homeomorphisms on the line. 

Moreover, if $\preceq$ is a relative order with respect to $C$, then there is an action $\rho_\preceq:G\to Homeo_+(\R)$ and a point $p\in \R$ such that $C=Stab_G(p)$ and
\begin{equation}\label{eq reference point} C\prec g C \Leftrightarrow p<\rho_\preceq(g)(p).\end{equation}}
\end{prop}

As with total orders, we call $\rho_\preceq$ a {\em dynamical realization} of the relative order $\preceq$ and call $p$ the {\em reference point} for the action.

\proof The proof is an adaptation of the  proof of \cite[Theorem 6.8]{ghys}, so we only give a sketch.
Suppose $G$ is acting non-trivially by orientation-preserving homeomorphisms of the line. There is some $g_0\in G$ and $x_0\in \mathbb{R}$ such that $g_0(x_0)\neq x_0$. Then $P:=\{g\in G\mid g(x_0)>x_0\}$ is a relative cone with respect to $C=Stab_G(x_0)\neq G$.

Conversely, suppose $\preceq$ is a $G$-invariant order on $G/C$. Then, since $G$ is countable, we can enumerate the cosets of $C$, say $g_0C,g_1C,\ldots$, and embed $(G/C,\preceq)$ into $(\R,\leq)$ via an order preserving map $t:G/C\to \R$.  The group $G$ then acts on $t(G/C)$, and, if $t$ is taken with a bit of care, as in \cite{ghys}, then this partial action can be extended to an action of $G$ by homeomorphisms of the line. Certainly, if $C=g_iC$, then $p:=t(g_iC)$ satisfies (\ref{eq reference point}).
 $\hfill\square$

\vsp

For each relatively convex subgroup $C$ one can consider  the set of all relative orders with respect to $C$, which we denote by $\mathcal{O}_C(G)$. More interestingly, one can consider the set of all relative orders on $G$: $\mathcal{O}_{Rel}(G)=\bigcup_C \mathcal{O}_C(G) $. Observe that $\mathcal{LO}(G)$, the set of left-orders on $G$, is contained in $\mathcal{O}_{Rel}(G)$. A natural topology can the defined on  $\mathcal{O}_{Rel}(G)$ by considering $\preceq $ as an element $\phi_{(\preceq,C)}\in \{\pm 1,*\}^G$ defined by 
$$\phi_{(\preceq,C)}: g\mapsto \left\{ \begin{array}{rl} 1 & \text{if } C\prec gC \\ 
-1 & \text{if } gC\prec C \\ * & \text{if }gC=C. \end{array}\right.$$
If we endow $\{\pm 1, *\}^G$ with the product topology, then $\{\pm 1,*\}^G$ is compact by Tychonoff's theorem and a basis of open neighborhoods of $\phi\in \{\pm 1,*\}^G$ consists of the sets of the form $U_{g_1,\ldots, g_n}(\phi):=\{\phi' \mid \phi'(g_i)=\phi(g_i)\}$, where $g_1,\ldots, g_n$ runs over all finite subsets of $G$. Note that $\{\pm 1,*\}^G$ is metrizable whenever $G$ is countable (see for instance \cite[\S 1]{navas}). For example if $G$ is finitely generated by $S$, then we can declare $dist(\phi_{(\preceq,C)},\phi_{(\preceq', C')})=1/2^n$ where $n$ is the largest integer such that  $\phi_{(\preceq,C)}$ and $\phi_{(\preceq', C')} $ agree on the ball $S^n$. Thus, in order to show that $\mathcal{O}_{Rel}(G)$ is compact we only need to show 

\begin{prop}\label{prop closed} {\em  If $G$ is finitely generated, then $\mathcal{O}_{Rel}(G)$ is closed inside $\{\pm 1,*\}^G$. In particular, $\mathcal{O}_{Rel}(G)$ is compact.}

\end{prop}

Before providing the proof, we exhibit an example showing that the finite generation hypothesis cannot be dropped.


\begin{exam} Consider $G=\oplus_{i\in \N} \Z$, and for every $n\in \N$ take $C_n=\oplus_{1\leq i\leq n} \Z$ and $\preceq_n$ a $G$-invariant order on $G/C_n$. Then $\phi_{(\preceq_n,C_n)}$ converges as $n\to \infty$ to the constant function $g\mapsto *$, which is certainly not a relative order of $G$ (since it is trivial).
\end{exam}
 
 

\noindent{\em Proof of Proposition \ref{prop closed}.}
Since $\{\pm 1,*\}^G$ is compact and metrizable, we only need to show that $\mathcal{O}_{Rel}(G)$ is sequentially closed. Let $\phi_n=\phi_{(\preceq_n, C_n)}$ be a convergent sequence of elements of  $\mathcal{O}_{Rel}(G)$, i.e. $\phi_n(g)$  converges  for all $g\in G$ (so, $\phi_n(g)$ is  an eventually constant sequence). Define $P:=\{g\in G\mid \phi_n(g)\to 1\}$ and $C:=\{g\in G\mid \phi_n(g)\to *\}$. Since $G$ is finitely generated, say by $S$, and $C_n=\phi_n^{-1}(*)$ are proper subgroups, there is a generator $g_0\in S$ such that $g_0\notin C_n$ for all sufficiently large $n$, and hence $g_0\notin C$. In particular, $C$ is a proper subgroup of $G$. Using Lemma \ref{lem:charconvex}, it is straightforward to check that $P$ is a relative cone with respect to $C$. $\hfill\square$

Let $P$ be a relative cone  with respect to $C$, then, for every $g\in G$,  $gPg^{-1}$ is a relative cone with respect to $gCg^{-1}$. More generally, for $g\in G$ and $\phi\in \{\pm1,*\}^G$, the map $g:\phi\mapsto \phi^g$, where \begin{equation} \label{eq conjugation}\phi^g(h)=\phi(ghg^{-1}),\end{equation}
sends bijectively the basic neighborhood $U_{g_1,\ldots,g_n}(\phi)$ onto $U_{g^{-1}g_1g,\ldots, g^{-1}g_ng}(\phi^g)$. Thus we have proved

\begin{prop} \label{prop action}{\em The natural conjugation action of $G$ on $\mathcal{O}_{Rel}(G)$ is an action by homeomorphisms.}

\end{prop}

The final assertion in Theorem \ref{teo Rel} is given by the next

\begin{cor}\label{coro fijos} {\em 
If $C$ is a relatively convex subgroup of $G$, then $C$ fixes a point in $\mathcal{O}_{Rel}(G)$. In fact $C$ fixes each point of $\mathcal{O}_C(G)$.}
\end{cor}

\proof Let $P$  be a relative cone with respect to $C$, and take $c\in C$. Then, by Proposition \ref{prop action} and Lemma \ref{lem:charconvex},  $G$ is the disjoint union of $cP^{-1}c^{-1} $, $C$, and $ cPc^{-1}$. But also, by Lemma \ref{lem:charconvex}, we have that $cPc^{-1}\subseteq P$, which implies that $cPc^{-1}=P$.
 $\hfill\square$ 

\vs

We finish this section with the example announced in the Introduction.

\begin{exam}\label{ejemplo} Let $H$ be a finitely generated, left-orderable group without homomorphisms onto $\Z$. For instance we can take $H$ being the perfect group with presentation $\langle a,b,c\mid a^2=b^3=c^7=abc\rangle$ (see \cite{Bergman, Thurston}). Let $G=\Z\times H$. Certainly $H$ is a relatively convex subgroup of $G$. We claim that $H$ does not preserve a probability measure when acting on $\mathcal{LO}(G)$. 

Indeed, since $\Z$ is amenable and the $\Z$-factor commutes with $H$, we conclude that if there is an $H$-invariant probability measure on $\mathcal{LO}(G)$, then an averaging procedure returns a $G$-invariant probability measure on $\mathcal{LO}(G)$. Following Morris' \cite[Remark 2.2]{Mo}, this is enough to show that $G$ admits a left-order of Conradian type $\preceq$. Since this condition clearly passes to subgroups, we conclude that $H$ admits a left-order of Conradian type, and therefore  admits  a surjective homomorphism onto $\Z$, which contradicts the choice of $H$.
\end{exam}

\section{Co-amenable relatively convex subgroup}
\label{sec coamenable}



In this section we prove Theorem \ref{teo main}. By Proposition \ref{prop real dyn} it is enough to show

\begin{thm}\label{teo witte} 
{\em Suppose $G$ is finitely generated and $C\subset G$ is a proper relatively convex subgroup which is also co-amenable. Then $G$ surjects onto $\Z$.}
\end{thm}

\subsection{Actions without crossings}

In addition to our work from Section \ref{sec relative orders}, our main tool is the concept of {\em crossings} for a group acting on a totally ordered space. This notion is due to  Beklaryan \cite{beklaryan} for the case of group action on the line, but here we use the analogous version for group actions on totally ordered space from \cite{Na-Ri}.

\begin{defi}\label{def crossing} Let $G$ be a group acting by order-preserving bijections of a totally ordered space $(\Omega,\leq)$. A {\em crossing} for the action is a 5-uple $(f,g;u,v,w)$, where $f,g\in G$, and $u,v,w\in \Omega$ such that:
\begin{itemize}
\item  $u < w < v$,

\item  $g^n (u) < v$  and  $u< f^n (v) $
 for every $n \in \mathbb{N}$,

\item  there exist $M,N$ in $\mathbb{N}$ so that 
$f^N (v) < w < g^M (u)$.

\end{itemize}

In a picture

\begin{center}
\begin{tikzpicture}

\draw[<->] (-3.2,-3.2)-> (3.2,3.2);
\draw[-] (-3,-3)-> (-3,3)->(3,3)->(3,-3)->(-3,-3);

\draw[fill] (-2.5,-3) circle (1.5 pt);
\draw[fill] (0,-3) circle (1.5 pt);
\draw[fill] (2.5,-3) circle (1.5 pt);

\draw (-2.5,-3.5) node{$u$};
\draw (0,-3.5) node{$w$};
\draw (2.5,-3.5) node{$v$};

\draw [-] (-2,-2) .. controls (-0,-1) .. (3,0); 
\draw (0,-1.7) node{$f$};

\draw [-] (-3,0) .. controls (0,1.5) .. (2,2); 
\draw (0,1.8) node{$g$};

\draw (0,-4.7) node{{\small Figure 1: The crossing $(f,g;u,v,w)$.}};

\end{tikzpicture}
\end{center}
\end{defi}

Let $G$ be a (perhaps not finitely generated)  group acting by order preserving bijections of a totally ordered space $(\Omega,\leq)$. For every $x\in \Omega$ and every $g\in G$, let $I_g(x)$ denote the convex envelope of the $\langle g\rangle$-orbit of $x$, that is $I_g(x)=\{y\in \Omega \mid g^n(x)\leq y \leq g^m(x)$ for some $n,m\in\Z\}$. 
Note that $I_g(x)=I_{g^{k}}(x)$ for every $k\in \mathbb{Z}-\{0\}$, $I_g(x)=I_g(y)$ for every $y\in I_g(x)$, and, for every $f,g\in G$, $f(I_g(x))=I_{fgf^{-1}}(f(x))$. It follows that 
\begin{equation}\label{eq I_g}
\text{if $I_g(x)=\Omega$ for some $x\in \Omega$, then $I_g(x)=\Omega$ for all $x\in \Omega$.} 
\end{equation}

Actions without crossings are sometimes referred to action by levels \cite{Na-Ri} or {\em emb\^oit\'es} \cite{triestino}. The reason for this is the following lemma whose proof can be found between the lines of \cite[Proof of Proposition 1.12]{Na-Ri}. For the reader's convenience we redo the argument here.  

\begin{lem}\label{lem emboite} {\em Suppose $G$ is a group acting by order preserving bijections  and without crossings on $(\Omega,\leq)$. Then either $I_g(x)$ and $I_f(y)$ are  disjoint or one is contained in the other.  }
\end{lem}

\proof Assume there are non-disjoint sets $I_g(x)$ and $I_f(y)$ neither of which contains the other. In particular, neither $g$ fixes $x$ nor $f$ fixes $y$. Without loss of generality  we may assume that $I_g(x)$ contains a point that is on the left of $I_f(y)$ (if this is not the case, just interchange the roles of $f$ and $g$) and therefore $I_f(y)$ contains a point on the right of $I_g(x)$. Moreover, by perhaps replacing $f$ and/or $g$ by their inverses, we may assume that $x'\leq g(x')$ for all $x'\in I_g(x)$ and $f(y')\leq y'$ for all $y'\in I_f(y)$. Take $u\in I_g(x)\setminus I_f(y)$, $w\in I_g(x)\cap I_f(y)$ and $v\in I_f(y)\setminus I_g(x)$. Then, one easily verifies that $(f,g;u,v,w)$ is a crossing for the action.  $\hfill\square$ 

\begin{prop}\label{prop fIg}
{\em Suppose $G$ is a group acting by order preserving bijections  and without crossings on $(\Omega,\leq)$. Then for every $f,g\in G$ and $x\in \Omega$, $f(I_g(x))$ either coincides with $I_g(x)$ or is disjoint from it. }
\end{prop}

\proof

Suppose that $I_g(x)\cap f(I_g(x))\neq \emptyset$. 
Let $g^f=fgf^{-1}$. 
Since $f(I_g(x))=I_{g^{f}}(f(x))$, by Lemma \ref{lem emboite} and replacing perhaps $f$ by $f^{-1}$, we can assume that $f(I_g(x))\subseteq I_g(x)$.

Suppose there is $y\notin I_g(x)$ such that $f(y)\in I_g(x)$, we will see that this implies that there is a crossing.
We will assume that $y>I_g(x)$, the case $y<I_g(x)$ is analogous. 
 Note that by Lemma \ref{lem emboite}, we have that $I_f(y)\supseteq I_g(x)$ and for all $z\in I_g(x)$, $f(z)<z$.
See Figure 2.

Observe that for all $n\in \mathbb{N}$,  $g^n(f(y))<f^{-1}(f(y))$ and thus $fg^n (f(y))<f(y)$. 
By taking $n$ large enough, we can assume that $fg^n (f^3(y))>f^2(y)$.
We claim that $I_{fg^n}(f(y))\cap f(I_g(x))\neq\emptyset$. Indeed, if the intersection was empty, then since $f^2(y)<f(y)$ we must have that  for all $z\in I_g(x)$,  $f(z)< I_{fg^n}(f(y))$, and since $f(I_g(x))=I_{(g^f)^n}(f(x))$ then $fg^n(z)<I_{fg^n}(f(y))$ for all $z\in I_g(x)$, which is a contradiction.

We let $v\in I_{fg^n}(f(y))\cap f(I_g(x))$. 
We have that $f(g^f)^n(v)< f^2(y)$ and by increasing $n$ if necessary, we can assume that $f(g^f)^n(f^3(y))>f^3(y)$.

Now it follows that $(fg^n, f(g^f)^n, f^3(y), f^2(y), v)$ is a crossing.

Therefore we get that for all $y\notin I_g(x)$, $f(y)\notin I_g(x)$ and since $f$ is a bijection,  we have that $f(I_g(x))=I_g(x)$.

\begin{center}
\begin{tikzpicture}

\draw[<->] (-3,-3)-> (3,3);

\draw (-2,1.3) node{$g$};
\draw [-] (-2.5,-2.5) .. controls (-3,1) .. (2.5,2.5);
\draw (2.2,3.3) node{$f^{-1}$};
\draw [-] (-2.5,-2.5) .. controls (-1,0) .. (2.3,4);
\draw (-1.5,0.3) node{$g^f$};
\draw [-] (-2.5,-2.5) .. controls (-2,1) .. (1.6,1.6);

\draw (0.3,0) node{$u$};
\draw[fill] (0.2,0.2) circle (1.5 pt);
\draw (1.5,1.1) node{$w$};
\draw[fill] (1.3,1.3) circle (1.5 pt);
\draw (2.5,2) node{$v$};
\draw[fill] (2.3,2.3) circle (1.5 pt);

\draw (0,-3.7) node{{\small Figure 2: $f$ {\em weakly contracts} $I_g$.}};

\end{tikzpicture}
\end{center}

$\hfill\square$

\begin{prop}\label{prop crossings} {\em Let  $G$ be a finitely generated group acting non-trivially by order-preserving bijections of a totally ordered space $(\Omega,\leq)$. Suppose that the action has no crossings. Then $G$ surjects onto $\Z$.}

\end{prop}

\proof Clearly, it is enough to restrict ourself to the case that the $G$-action on $\Omega$ has all orbits unbounded  in the sense that for all $x< y$ there is $g\in G$ such that $y< g(x)$ (in the case $\Omega=\R$ this is just assuming that the action has no global fixed points). Assuming this, we will say that $g\in G$ is {\em cofinal} if for every $x< y$, there is $n\in \Z$ such that $y< g^n(x)$ (in the case $\Omega=\R$, this is saying that $g $ has no fixed points). It follows from \eqref{eq I_g} that $g$ is cofinal if and only if $I_g(x)=\Omega$ for all $x\in \Omega$. 

Let $F_G$ denote the set of elements in $G$ which are not cofinal. This is, in general, a subset invariant under conjugation, but, since the $G$-action has no crossings, we have that $F_G$ is a normal subgroup (see  \cite[Proposition 2.13]{Na-Ri}). 

We claim that  $F_G$ is a proper subgroup. To prove the claim observe that $G$ is finitely generated, so we let $S=S^{-1}$ be a finite generating set of $G$. We also fix $x_0\in \Omega$ and let $g_0\in S$ be such that $g(x_0)\leq g_0(x_0)$ for all $g\in S$. Thus $\emptyset\neq I_g(x_0)\cap I_{g_0}(x_0)$ for every $g\in S$, and it follows from Proposition \ref{prop fIg} that $I_{g_0}(x_0)$ is invariant under every $g\in S$ and therefore it is invariant under $G$. Since the $G$-action is unbounded,  $I_{g_0}(x_0)=\Omega$ and the claim follows.

We now claim that $G/F_G$ is torsion-free abelian. For this we let $I^*$ be the union of all the intervals $I_f(x_0)$ for $f\in F_G$. By Lemma \ref{lem emboite}, we have that the union defining $I^*$ is nested and so $I^*$ is an interval. Let $f\in F_G$ and $g\in G\setminus F_G$. Suppose that $g(I_f(x_0))=I_f(x_0)$. Then $I_g(x_0)=I_f(x_0)\neq \Omega$, a contradiction. From Proposition \ref{prop fIg}, we deduce that $g(I_f(x_0))$ is disjoint from $I_f(x_0)$. Since $I^*$ is a nested union of intervals of the form $I_f(x_0)$, $f\in F_G$, we conclude that $g(I^*)$ is disjoint from $I^*$ for every $g\notin F_G$ and that $f(I^*)=I^*$ for every $f\in F_G$. Finally, since $F_G$ is a normal subgroup, we have that $G/F_G$ acts on the $G$-orbit of $I^*$, and the kernel of this new action is trivial. We can therefore {\em induce} a (total) left-order on $G/F_G$ by declaring that $\overline{g}\succ \overline{id}$ if and only if $g(I^*)$ is to the right of $I^*$. This is an Archimedean left-order on $G/F_G$, meaning that for every $\overline{f},\overline{g}\in G/F_G$ we have that there is $n\in \Z$ such that $\overline{f}\prec \overline{g}^n$. It then follows from Holder's theorem (see \cite[\S 3.1]{GOD}) that $G/F_G$ is a subgroup of $(\R,+)$, which implies the claim.

Since $G$ is finitely generated, $G/F_G$ is a finitely generated torsion-free abelian group and hence maps onto $\mathbb{Z}$.
$\hfill\square$

\subsection{Right recurrent orders}
The next definition is due to Morris for the case of left-orders \cite{Mo}.

\begin{defi} 
Let $\preceq$ be a relative order of $G$ with respect to $C$. We say that $\preceq$ is {\em recurrent for every cyclic subgroup} if for every $h\in G$ and every finite sequence $\lambda_1,\dots, \lambda_k$ in $G$ such that  $\lambda_1C\prec \ldots \prec \lambda_kC$, there exist positive integers $n_i\to \infty$ such that 
$$\lambda_1  h^{n_i} C\prec \ldots \prec \lambda_k h^{n_i}C.$$
\end{defi}

\begin{lem}\label{lem conrad} 
{\em Suppose $\preceq$ is a relative order of $G$ with respect to $C$ that is recurrent for every cyclic subgroup. Then the left-action of $G$ on $(G/C,\preceq)$ is an action without crossings. In particular, if $G$ is finitely generated, then $G$ surjects onto $\Z$. }
\end{lem}

\proof  Suppose that the $G$-action on $\Omega=G/C$ has a crossing, say $(f,g;uC,vC,wC)$. In particular we have that
$$f^Ng^MwC\prec wC.$$
Let $h=w^{-1}f w$, and let $k>N$. Then by the third condition in the definition of crossing we have that $wC\prec g^Mf^k w C=g^M w h^k C $. But then, by $G$-invariance we have that for $k>N$, 
$$f^Ng^Mwh^k C =f^Ng^Mf^k w C \succ f^N wC \succ f^kwC=wh^kC,$$
contradicting the right-recurrence for $\langle h \rangle$. This shows the first assertion of the lemma. The second one follows from Proposition \ref{prop crossings} by taking $(\Omega,\leq)=(G/C,\preceq)$. $\hfill\square$

\vsp

We can now give the

\vsp

\noindent {\bf Proof of Theorem \ref{teo witte}:} Let $G$ and $C$ be as in the statement, and consider the conjugation action of $G$ on $\mathcal{O}_{Rel}(G)$ denoted by $g:\phi_{(\preceq,C)}\mapsto \phi^g_{(\preceq,C)}$ (see equation (\ref{eq conjugation})). By Theorem  \ref{teo Rel} the space $\mathcal{O}_{Rel}(G)$ is compact and $C$ acts with a fixed point. In particular, since $C$ is co-amenable, there is a $G$-invariant probability measure $\mu$ on $\mathcal{O}_{Rel}(G)$.

Now, the Poincar\'e Recurrence Theorem (see, for instance, \cite{sinai}) implies that for every $\phi_0=\phi_{(\preceq,C)}$ in the support of $\mu$, every neighborhood $U_{f_1,\ldots, f_n}(\phi_0)$ of $\phi_0$ and every $g\in G$, there is  a $\mu$-null set $Z_{g,f_1,\ldots,f_n}$ such that for every $\phi\in U_{f_1,\ldots,f_n}(\phi_0)\setminus Z _{g,f_1,\ldots,f_n}(\phi)$ there is a sequence $n_i\to \infty$ such that $\phi^{g^{-n_i}}$ belongs to $U_{f_1,\ldots,f_n}(\phi_0)$.  But since $G$ is countable and $\mathcal{O}_{Rel}(G)$ is metrizable (hence it has a countable basis of open neighborhoods) and compact (hence it has a dense countable subset), we can take a countable union of the above $Z$'s to conclude:
\begin{itemize}
\item[(Recurrence)] {\em There is $Z\subset \mathcal{O}_{Rel}(G)$ with $\mu(Z)=0$ such that for every $\phi=\phi_{(\preceq,C)} \in \mathcal O_{Rel}(G)\setminus Z$, every open neighborhood $U_{f_1,\ldots,f_n}(\phi)$ of $\phi$ and every $g\in G$, there is a sequence $n_i\to\infty$ such that
$\phi^{g^{-n_i}}\in U_{f_1,\ldots,f_n}(\phi).$}
\end{itemize}

Let $\phi_{(\preceq,C)}$ be a relative order in $\mathcal{O}_{Rel}(G)\setminus Z$. We claim that $\preceq$ is right  recurrent for every cyclic subgroup. Indeed, let $\lambda_1C \prec \ldots \prec \lambda_nC$,  $g\in G$, and
consider the neighborhood $U=U_{\lambda_1^{-1}\lambda_2, \lambda_2^{-1}\lambda_3, \ldots,\lambda_{n-1}^{-1}\lambda_n}(\phi)$ of $\phi$. Note that $U$ coincides with the set made of all relative orders on which $\lambda_k$ is {\em strictly smaller} than $\lambda_{k+1}$ for every $k=1,\ldots,n-1$. Now, by (Recurrence), we have that there is $n_i\to \infty$ such that $\phi^{g^{-n_i}}\in U$. This means that $\phi^{g^{-n_i}}(\lambda_k^{-1}\lambda_{k+1})= \phi(\lambda_k^{-1}\lambda_{k+1})=1$. But by definition $\phi^{g^{-n_i}}(\lambda_k^{-1}\lambda_{k+1})= \phi(g^{-n_i}\lambda_k^{-1} \lambda_{k+1} g^{n_i})$, thus $1=\phi(g^{-n_i}\lambda_k^{-1} \lambda_{k+1} g^{n_i})$ which implies $C \prec g^{-n_i}\lambda_k^{-1} \lambda_{k+1} g^{n_i}C$. Using left-multiplication invariance and transitivity of $\preceq$ we conclude that $\lambda_1g^{n_i}C\prec \ldots \prec \lambda_ng^{n_i}C$, so $\preceq$ is right recurrent for every cyclic subgroup.  In particular, Lemma \ref{lem conrad} implies that $G$ surjects  onto $\Z$. $\hfill\square$

\end{document}